\documentclass[10pt]{article}

\usepackage{amssymb}
\usepackage{amsfonts}
\usepackage{latexsym}
\usepackage[T1,T2A]{fontenc}
\usepackage[cp1251]{inputenc}
\usepackage[T2B]{fontenc}
\usepackage[bulgarian,russian,english]{babel}

\usepackage{color}

\newtheorem {problem}{\bf Problem}

\newtheorem{proposition}{\bf Proposition}

\title {Bitwise operations related to a combinatorial problem on binary matrices
 }
\author {Krasimir Yankov Yordzhev}
\date {}

\begin {document}
\inputencoding{cp1251}

\maketitle
\begin{abstract}
Some techniques for the use of bitwise operations are described in the article.
As an example, an open problem of isomorphism-free generations of combinatorial objects is discussed.
An equivalence relation on the set of square binary matrices having the same number of units in each row and each column is defined.
Each binary matrix is represented using ordered n-tuples of natural numbers.
It is shown how by using the bitwise operations can be implemented an algorithm that gets  canonical representatives which are extremal elements of equivalence classes  relative  to  a double  order on the set of considered objects.
\end{abstract}

\section{Introduction}\label{intr}

The present study is thus especially useful for students educated to become programmers as well as for their lecturers. A  meaningful example for the advantages of using bitwise operations for creating effective algorithms in programming is presented in this article. We will consider an open combinatorial problem on binary matrices and its solution using the algorithm for some values of the integer parameters $n$ and $k$. To implement the algorithm, we will use essentially bitwise operations.

The use of bitwise operations is a powerful method used in C/C++ and Java programming languages. Unfortunately, in the widespread books on this topic there is incomplete or no description for the work of the bitwise operations. The aim of this article is to correct this lapse to a certain extent and present a meaningful example of a programming task, where the use of bitwise operations is appropriate in order to facilitate the work and to increase the effectiveness of the respective algorithm.

This work is an extension and complement to  \cite{umb2009}.

A \emph{binary} (or   \emph{boolean}, or (0,1)-\emph{matrix}) is a matrix whose all elements belong to the set $\mathcal{B} =\{ 0,1 \}$. With $\mathcal{B}_n$ we will denote the set of all  $n\times n$  binary matrices.

Some algorithms for isomorphism-free generations of combinatorial objects are discussed in detail in \cite{Bouyukliev}.
In our work we will consider a problem of this type. Its formulation is as follows: A set of binary matrices $\mathcal{L}\subseteq \mathcal{B}_n$ is given. In  $\mathcal{L}$ is defined an equivalence relation.  An algorithm which did not study every element of the set $\mathcal{L}$, and which receives one representative of each equivalence class  to be described. For this purpose, we will use significantly bitwise operations.

In Section \ref{prblms} we formulate the problem and we give some well known results. We will describe in detail an algorithm for computer solution of the formulated problem in Section \ref{sect3}. Section \ref{bitwiseop} is only for reference.

\section{Preliminaries and problem formulation}\label{prblms}

Let $n$ and $k$ be positive integers. We let $\Lambda_n^k$ denote the set of all $n\times n$ binary matrices in each row and each column of which there are exactly $k$ in number 1's.
Let us denote with $\lambda (n,k) =|\Lambda_n^k |$ the number of all elements of $\Lambda_n^k$.

There is not any known formula to calculate  the $\lambda (n,k)$  for all $n$ and $k$. There are formulas for the calculation of the function $\lambda (n, k)$ for each $n$ for relatively small values of $k$; more specifically, for $k = 1$, $k = 2$ and $k = 3$. We do not know any formula to calculate the function $\lambda (n, k)$  for $k > 3$ and for all positive integer $n$.

It is easy to prove the following well-known formula

$$
\lambda (n,1) =n!
$$

The following formula

$\displaystyle
\lambda (n,2) =$ $\displaystyle \sum_{2x_2 +3x_3 + \cdots +nx_n =n}$ $\frac{(n!)^2}{\displaystyle \prod_{r=2}^n x_r !(2r)^{x_r}}$

is well known \cite{6}.

One of the first recursive formulas for the calculation of
$\lambda (n,2)$ appeared in \cite{anand} (see also \cite[p. 763]{gupta}).

$\displaystyle \lambda (n,2) =
\frac{1}{2} n(n-1)^2 \left[ (2n-3) \lambda (n-2,2) +(n-2)^2 \lambda (n-3,2) \right]$ for  $n \ge 4$

$$\displaystyle \lambda (1,2) =0,\quad \lambda (2,2) =1,\quad
\lambda (3,2) =6$$

Another recursive formula for the calculation of $\lambda (n,2)$ occurs  in \cite{13}.

$\displaystyle \lambda  (n,2) = (n-1)n\lambda (n-1,2)
+\frac{(n-1)^2 n}{2} \lambda (n-2,2)$ for  $n\ge 3$

$\lambda (1,2) =0$, $\lambda (2,2) =1$

The next recursive system is  to calculate $\lambda (n,2)$.

$\lambda (n+1,2) =$

$n(2n-1)\lambda (n,2) + n^2 \lambda (n-1,2) - \pi(n+1);$ $n\ge 2$

$\pi (n+1) =$

${n^2 (n-1)^2 \over 4} [8(n-2)(n-3) \lambda (n-2,2) +(n-2)^2 \lambda (n-3,2) -4\pi (n-1)] ;$  $ n\ge 4$

$\lambda (1,2) = 0 $, $\lambda (2,2) = 1 $, $\pi (1) = \pi (2) = \pi (3) = 0 $, $\pi (4) = 9$

where  $\pi (n)$ identifies the number of a special class of $\Lambda_n^2$-matrices \cite{iord}.

The following formula is an explicit form for the calculation of $\lambda (n,3)$.

$$
 \lambda (n,3)=\frac{n!^2}{6^n} \sum \frac{(-1)^\beta (\beta
+3\gamma )! 2^\alpha 3^\beta}{\alpha !\beta ! \gamma !^2 6^\gamma}
$$
where the sum is done as regards all $ \frac{(n+2)(n+1)}{2}$ solutions in nonnegative integers of the equation   $\alpha +\beta+\gamma =n$ \cite{stein}. As it is noted in  \cite{stan}, the above formula does not give us good opportunities to study behavior of  $\lambda (n,3)$.

Let $A,B\in\Lambda_n^k$. We will say that $A\sim B $, if $A$ is obtained from $B$ by moving some rows and/or columns. Obviously,  the relation  defined like that is an equivalence relation. We denote with
$$\mu (n,k) = \left| {\Lambda_n^k}_{/_\sim} \right|$$
the number of equivalence classes on the above defined relation.

\begin{problem}
Find $\mu(n, k)$ for given integers $n$ and $k$, $1\le k<n$.
\end{problem}

The task of finding the number of equivalence classes for all integers $n$ and $k$, $1\le k<n$ is an open scientific problem.
We partially solve this problem by making a computer program to find this number for some (not great) values of $n$ and $k$.
Moreover, using bitwise operations, our algorithm will receive one representative from each equivalence class without examining the whole set
$\Lambda_n^k$.

\section{Bitwise operations}\label{bitwiseop}
Bitwise operations can be applied for integer data type only, i.e. they cannot be used for float and double types. For the definition of the bitwise operations and some of their elementary applications could be seen, for example, in \cite{Davis,Kernigan,Schildt}.

We assume, as usual that bits numbering in variables starts from right to left, and that the number of the very right one is 0.

Let  \verb"x,y"  and \verb"z" are integer variables or constants of one type, for which   bits are needed. Let \verb"x"  and  \verb"y" are initialized (if they are variables) and let the  assignment \verb"z = x & y;" (\emph{bitwise AND}), or \verb"z = x | y;" (\emph{bitwise inclusive OR}), or \verb"z = x ^ y;" (\emph{bitwise exclusive OR}), or \verb"z = ~x;" (\emph{bitwise NOT}) be made. For each $i=0,1,2,\ldots ,w-1$,  the new contents of the $i$-th  bit in \verb"z"  will be as it is presented in the Table \ref{bitwise}.

\begin{table} [h]
\begin{center}
\begin{tabular}{||c|c|c|c|c|c||}
  \hline\hline
 $i$-th bit of  &  $i$-th bit of  &  $i$-th bit of  & $i$-th bit of  &  $i$-th bit of  &  $i$-th bit of  \\
 \verb"x" & \verb"y" & \verb"z = x & y;" & \verb"z = x | y;" & \verb"z = x ^ y;" & \verb"z = ~x;"\\
\hline\hline
  0 & 0 & 0 & 0 & 0 & 1\\ \hline
  0 & 1 & 0 & 1 & 1 & 1\\ \hline
  1 & 0 & 0 & 1 & 1 & 0\\ \hline
  1 & 1 & 1 & 1 & 0 & 0\\ \hline
  \hline
\end{tabular}
\caption{Bitwise operations}\label{bitwise}
\end{center}
\end{table}

In case that  \verb"k"  is a nonnegative integer, then the statement \verb"z = x<<k"  (\emph{bitwise shift left}) will write $(i+k)$ in the bit of \verb"z" the value of the $k$ bit of \verb"x", where $i=0,1,\ldots ,w-k-1$, and the very right $k$  bits of \verb"x" will be filled by zeroes. This operation is equivalent to a multiplication of  \verb"x"  by  $2^k$.

The statement \verb"z=x>>k" (\emph{bitwise shift right}) works the similar way. But we must be careful if we use the programming language C or C + +, as in various programming environments this operation has different interpretations – somewhere $k$  bits of \verb"z"  from the very left place are compulsory filled by 0 (logical displacement), and elsewhere the very left $k$  bits of \verb"z"  are filled with the value from the very left (sign) bit; i.e. if the number is negative, then the filling will be with 1 (arithmetic displacement). Therefore it is recommended to use \verb"unsigned" type of variables (if the opposite is not necessary) while working with bitwise operations. In the Java programming language, this problem is solved by introducing the two different operators: \verb"z=x>>k" and \verb"z=x>>>k" \cite{Schildt}.

Bitwise operations are left associative.

The priority of operations in descending order is as follows:  \verb"~" (bitwise NOT); the arithmetic operations * (multiply), / (divide), \% (remainder or modulus); the arithmetic operations + (addition) - (subtraction); the  bitwise operations $<<$ and $>>$; the relational  operations $<$, $>$, $<=$, $>=$, ==, !=; the bitwise operations \&,\verb"^" and $|$; the logical operations \&\& and ||.

To compute the value of the $i$-th  bit of an integer variable \verb"x" we can use the function:
\begin{verbatim}
int BitValue(int x, unsigned int i) {
	return ( (x & 1<<i) == 0 ) ? 0 : 1;
}
\end{verbatim}

The next function prints an integer in binary notation. We don't consider and we don't print the sign of integer. For this reason we work with $|n|$

\begin{verbatim}
void DecToBin(int n)
{
     n = abs(n);
     int b;
     int d = sizeof(int)*8 - 1;
     while ( d>0 && (n & 1<<(d-1) ) == 0 ) d--;
     while (d>=0)
           {
           b= 1<<(d-1) & n ? 1 : 0;
           cout<<b;
           d--;
           }
}
\end{verbatim}

The following function calculates the number of 1's in the binary representation of an integer n.
Again we ignore the sign of the number.

\begin{verbatim}
int NumbOf_1(int n)
{
     n = abs(n);
     int temp=0;
     int d = sizeof(int)*8 - 1;
     for (int i=0; i<d; i++)
         if (n & 1<<i) temp++;
     return temp;
}
\end{verbatim}

\section{Description and implementation of the algorithm}\label{sect3}
Let $\mathbb{N} $ be the set of natural numbers and let
$$\mathcal{T}_n =\left\{ \langle x_1 ,x_2 ,\ldots ,x_n \rangle \; |\; x_i \in \mathbb{N},\; i=1,2,\ldots ,n\right\}$$

An one to one corresponding
$$\varphi \; :\; \mathcal{B}_n \stackrel{\sim}{\longrightarrow} \mathcal{T}_n$$
which is based on the binary presentation of the natural numbers, is described in \cite{umb2009}. If $A\in \mathcal{B}_n$ and $\varphi (A) =\langle x_1 ,x_2 ,\ldots x_n \rangle $, then $i$-th row of  $A$ is integer $x_i$ written in binary notation.

In \cite{Kostadinova}, it is proved that the representation of the elements of $\mathcal{B}_n$ using ordered $n$-tuples of natural numbers leads to making a fast and saving memory algorithms.

Let $A\in \mathcal{B}_n$ and let $\mathbf{x} =\langle x_1 ,x_2 ,\ldots ,x_n \rangle =\varphi (A)$. Then we denote
$$\mathbf{x}^t =\varphi (A^t ) ,$$
where $A^t \in \mathcal{B}_n$ is the transpose of the matrix  $A$.

Let $\mathbf{x}= \langle x_1 ,x_2 , \ldots ,x_n \rangle$  and let $\mathbf{x}^t = \langle y_1 ,y_2 ,\ldots ,y_n \rangle $. $\bf x$ we will call {\it canonical element}, if $x_1\le x_2 \le \cdots \le x_n$  and $y_1\le 2_2 \le \cdots \le y_n$.

\begin{proposition}\label{prop1}
There is un unique canonical element in every equivalence class of factor-set  ${\Lambda_n^k}_{/_\sim}$.
\end{proposition}

The proof of proposition \ref{prop1} is within the reach of any student who has successfully studied the properties of the binary system concept and we will miss it here.

Proposition 1 is the base of our algorithm, which we describe in brief below. For its implementation, we will use also the functions shown in section 3.

As it is well known, there are exactly $2^n$ nonnegative integers, which are presented with no more than $n$ digits in binary notation.
We  need to select all of them, which have exactly $ k $ 1's in binary notation.
Their number is ${n\choose k}\ll 2^n$. We could use the function \verb"NumbOf_1(int)" from section  \ref{bitwiseop},
but then we have to use it for each integer from the interval  $[0,2^n -1]$, i.e. $2^n$ times.
We will describe an algorithm that directly receives the necessary elements without checking whether any integer $m\in [0,2^n -1]$
satisfies the conditions. We will remember the result in the array \verb"p[]" of size $c={n\choose k}$.
Moreover, the obtained array is sorted in ascending order and there are no duplicate elements.
The algorithm is based on the fact that the set of all ordered $m$-tuples $\mathcal{B}^m =\langle b_1 ,b_2 ,\ldots , b_m \rangle  $,
$b_i \in \mathcal{B} =\{ 0,1\}$, $i=1,2,\ldots ,m$, $m=1,2,\ldots ,n$,  is partitioned into two disjoint subsets
$\mathcal{B}^m =\mathcal{M}_1 \cup \mathcal{M}_2$, $\mathcal{M}_1 \cap \mathcal{M}_2 =\emptyset$,
where $\mathcal{M}_1 =\{ \langle b_1 ,b_2 ,\ldots , b_m \rangle \; |\; b_1 =0\}$ and
$\mathcal{M}_2 =\{ \langle b_1 ,b_2 ,\ldots , b_m \rangle \; |\; b_1 =1\}$. The described recursive algorithm again uses bitwise operations.

\begin{verbatim}
void DataNumb(int p[], unsigned int n, unsigned int k, int& c)
{
     if (k==0)
        {
        c = 1;
        p[0] = 0;
        }
        else if (k==n)
          {
          c = 1;
          p[0] = 0;
          for (int i=0; i<k; i++) p[0] |=  1<<i;
          }
          else
               {
               int p1[10000], p2[10000];
               int c1, c2;
               DataNumb(p1, n-1, k, c1);
               DataNumb(p2, n-1, k-1, c2);
               c = c1+c2;
               for (int i=0; i<c1; i++) p[i] = p1[i];
               for (int i=0; i<c2; i++) p[c1+i] = p2[i] | 1<<(n-1);
               }
}
\end{verbatim}

We also will use bitwise operations in constructing the next two functions.

The function \verb"int n_tuple(int[], int, int, int)" gets all $ t = {n + k-1 \choose k} $ (combinations with repetitions)
ordered $n$-tuples  $\langle x_1, x_2, \dots , x_n \rangle$, where  $0 \le x_1 \le x_2 \le \ldots \le x_n <c$,  $x_i$, $i=1,2,\ldots ,n$
are elements of sorted array \verb"p[]" of size \verb"c". As a result, the function returns the number of canonical elements.

The function \verb"bool check(int[], int)" refers to the use of each received $n$-tuples.
It examines whether this is a canonical element and prints it.

\begin{verbatim}
bool check(int x[], int n, int k)
{
     int yj; // the integer representing column (n-j)
     int y0=0; // integer preceding column j
     int b;
     for (int j=n-1; j>=0; j--)
     {
         yj=0;
         for (int i=0; i<n; i++)
         {
             b = 1<<j & x[i] ? 1 : 0;
             yj |= b << (n-1-i);
         }

         if (yj<y0 || (NumbOf_1(yj) != k)) return false;
         y0 = yj;
     }
     // We have received a canonical element. Print it:
     for (int i=0; i<n; i++) cout<<x[i]<<"  ";
     cout<<'\n';
     return true;
}

int n_tuple(int p[], int n, int k, int c)
{
     int t=0;
     int a[n], x[n];
     int indx = n-1;
     for (int i=0; i<n; i++) a[i]=0;
     while (indx >= 0)
     {
           for (int i=indx+1; i<n; i++) a[i] = a[indx];

           for (int i=0; i<n; i++) x[i] = p[a[i]];
           if(check(x,n,k)) t++;

           indx = n-1;
           a[indx]++;
           while (indx>=0 && a[indx]==c)
           {
                 indx--;
                 a[indx]++;
           }
     }
     return t;
}
\end{verbatim}

The description of the main function, we leave to the reader.

\begin{table}[h]
\begin{center}
\begin{tabular}{||lr||c|c|c|c|c|c|c|c||}
  \hline\hline
       & $n$ & 2 & 3 & 4 & 5 & 6 & 7 & 8 & 9  \\
  $k$  &     &   &   &   &   &   &   &  &   \\
  \hline\hline
  1 &   & 1 & 1 & 1 & 1 & 1 & 1 & 1 & 1  \\
  \hline
  2 &   &   & 1 & 2 & 5 & 13 & 42 & 155 & 636 \\
  \hline
  3 &   &   &   & 1 & 3 & 25 & 272 & 4 070 & 79 221 \\
  \hline
  4 &   &   &   &   & 1 & 5 & 161 & 7 776 & 626 649 \\
  \hline
  5 &   &   &   &   &   & 1 & 8 & 1 112 & 287 311 \\
  \hline
  6 &   &   &   &   &   &   & 1 & 13 & 8 787 \\
  \hline
  7 &   &   &   &   &   &   &   & 1 & 21 \\
  \hline
  8 &   &   &   &   &   &   &   &   & 1 \\
  \hline\hline
\end{tabular}
\caption{The number of equivalence classes for $1\le k<n \le 9$}\label{concl}
\end{center}
\end{table}

\section{Conclusion}
The number of equivalence classes for $1\le k<n \le 9$ are given in Table \ref{concl}, which is obtained through the work of the algorithms described in this paper.

The ideas described in this article can be used for finding the cardinality of other factor-sets of binary matrices

Krasimir Yankov Yordzhev

South-West University  ''N. Rilsky''

Faculty of Mathematics and Natural Sciences

2700 Blagoevgrad

Bulgaria

\rm e-mail:  yordzhev@swu.bg
\end{document}